\documentclass{amsart}
\newtheorem{theorem}{Theorem}[section]

\theoremstyle{definition}

\newtheorem{example}[theorem]{Example}

\theoremstyle{remark}

\numberwithin{equation}{section}

\usepackage{graphicx}
\usepackage{amsfonts}
\usepackage{color}
\everymath{\displaystyle}

\begin{document}

%% Title
\title[An Experimental Mathematics Perspective]{An Experimental Mathematics Perspective on the Old, and still Open, Question of When To Stop?}

%% Author # 1
\author{Luis A. Medina}
\address{Department of Mathematics, University of Puerto Rico, San Juan, PR 00931}
\email{luis.medina@uprrp.edu}

%% Author # 2
\author{Doron Zeilberger}
\address{Department of Mathematics, Rutgers University, Piscataway, NJ 08854}
\email{zeilberg@math.rutgers.edu}

\maketitle

%% Note at the beginning.
\noindent
{\footnotesize``You got to know when to hold em, know when to fold em, 
know when to walk away... ''
\begin{flushleft}
\quad\quad -Kenny  Rogers
\end{flushleft}
}

\bigskip

\noindent
{\bf Supporting Software}:
This article is accompanied by Maple packages {\tt ChowRobbins}, {\tt STADJE}, and { \tt WALKSab},
and Mathematica packages {\tt Builder.m} (and notebook {\tt Builder.nb}) as well as {\tt STADJE.m},  available
from the webpage of this article

{\tt http://www.math.rutgers.edu/\~{}zeilberg/mamarim/mamarimhtml/stop.html}{}.

\bigskip
%%%%%%%%%%%%%%%%%%%%%%%%%%%%%%%%%%%%%%%%%%%%%%%%%%%%%%%%%%%%%%%%%%%%%%%%%%%%%%%%%%%%%%%%%%%%%%%%%
% Section: When to Stop?
%%%%%%%%%%%%%%%%%%%%%%%%%%%%%%%%%%%%%%%%%%%%%%%%%%%%%%%%%%%%%%%%%%%%%%%%%%%%%%%%%%%%%%%%%%%%%%%%%

\section{When to Stop? }

In a delightful and insightful recent ``general'' article \cite{hill1}, the great probabilist and master expositor Theodore Hill described, amongst numerous other intriguing things, a more than forty-year-old open problem, due to Y.H. Chow and Herbert Robbins \cite{chowrob} that goes as follows:\\

{\it Toss a fair coin repeatedly and stop whenever you want, receiving as a reward the average number of heads accrued at the time you stop. 
If your first toss is a head, and you stop, your reward is $1$ Krugerrand. Since you 
can never have more than $100$ percent heads, it is clearly optimal to stop in that case. 
If the first toss is a tail, on the other hand, it is clearly best not to stop, since your reward would be zero...}\\

\noindent
Then Ted Hill goes on to comment that if the first toss is a tail and the second is a head, then it is good to {\bf go}, since by the law of large numbers, you would eventually do (at least slightly) better than one half. [It turns out that in this case of one head and one tail, the expected gain of continuing the game is larger than $0.6181$].

Hill further claims that it is optimal to {\bf stop} if the initial sequence is tail-head-head. [This is {\it wrong}. It turns out, thanks to our computations, that it is optimal to {\bf go}, and the expected gain is $> 0.6693$ rather than $2/3$.]

The exact stopping rule, i.e. the decision whether to {\bf stop} or {\bf go}, is still an open problem for (infinitely) many cases. As we will see, it is easy (with computers!)  to prove that it is optimal to {\bf go} for many cases where this is indeed the case, but proving {\it rigorously} that 
for a given position it is optimal to {\bf stop} is a challenging, still open, problem. It is analogous to disproving vs. proving a mathematical conjecture. To disprove it, all you need is to come-up with a specific counterexample, whereas to prove it, you need to show that no counterexample exists.

%%%%%%%%%%%%%%%%%%%%%%%%%%%%%%%%%%%%%%%%%%%%%%%%%%%%%%%%%%%%%%%%%%%%%%%%%%%%%%%%%%%%%%%%%%%%%%%%%
% Section: The Continuous Limit
%%%%%%%%%%%%%%%%%%%%%%%%%%%%%%%%%%%%%%%%%%%%%%%%%%%%%%%%%%%%%%%%%%%%%%%%%%%%%%%%%%%%%%%%%%%%%%%%%
\section{The Continuous Limit}

Way back in the mid sixties, this problem  was tackled by such luminaries as Chow and Robbins themselves \cite{chowrob}, Aryeh Dvoretzky \cite{dvoretzky1}, and Larry Shepp \cite{shepp1}. Chow and Robbins proved the existence of a stopping sequence, $\beta_n$, such that you stop as soon the number of heads minus the number of tails, after $n$ tosses, is $\geq \beta_n$. While Chow and Robbins only proved the {\it existence} of the ``stopping sequence'', Dvoretsky \cite{dvoretzky1} proved that $\beta_n/\sqrt{n}$ lies between two constants, for $n$ sufficiently large, while Larry Shepp \cite{shepp1} went further and proved  that 
\begin{equation}
\lim_{n \rightarrow {\infty}}\frac{\beta_n}{\sqrt{n}}
\end{equation}
exists and equals $0.83992 \dots$, a root of a certain transcendental equation.

But this beautiful work, like most of ``modern'' probability theory,  is {\it asymptotic}, talking about large $n$. It tells us nothing, for example,  about the still open $\beta_8$ (presumably $2$)  and not even about $\beta_{100}$. For example, the still-open question whether $\beta_8=2$  can be phrased as follows.\\

\noindent
{\it If currently you have five heads and three tails, should you {\bf stop}? }\\

\noindent
If you {\bf stop}, you can {\it definitely} collect $5/8=0.625$, whereas if you keep {\bf go}ing, your expected gain is $>0.6235$, but no one currently knows to {\it prove} that it  would not eventually exceeds $5/8$ (even though this seems very unlikely, judging by numerical heuristics).

%%%%%%%%%%%%%%%%%%%%%%%%%%%%%%%%%%%%%%%%%%%%%%%%%%%%%%%%%%%%%%%%%%%%%%%%%%%%%%%%%%%%%%%%%%%%%%%%%
% Section: The Role of Computers in Pure Mathematics Research 
%%%%%%%%%%%%%%%%%%%%%%%%%%%%%%%%%%%%%%%%%%%%%%%%%%%%%%%%%%%%%%%%%%%%%%%%%%%%%%%%%%%%%%%%%%%%%%%%%
\section{ The Role of Computers in Pure Mathematical Research}

We really enjoyed Hill's fascinating article, but we beg to differ on one (important!) issue. Hill (\cite{hill1}, p. 131) claims that:\\

\noindent
{\it ``Computers were not useful for solving that problem. In fact, all the problems described in this article were solved using traditional mathematicians' tools-working example after example with paper and pencil; settling the case for two, three, and then four unknowns; looking for patterns; waiting for the necessary {\bf Aha!} insights; and then searching for formal proofs in each step.''}\\

So far, this is all {\it factual}, so there is nothing to disagree with. Ted Hill was merely describing how he and his colleagues do research in pure mathematics. But then came an {\it opinion} that we do {\bf not} agree with:\\

\noindent
{\it ``Computers are very helpful for after-the-fact applications of many results, such as backward induction. But in theoretical probability, computers often do not significantly aid the discovery process.''}\\

This may have been true in the {\it past}, and to a large extent still at {\it present}, but we believe that in the {\it future} computers will be more and more useful even-and perhaps especially-in theory, since in addition to their obvious role as {\bf number-crunchers}, they are also starting to
do a great job as {\bf symbol-crunchers}, and even as {\bf idea-crunchers}. One recent example is \cite{zeilberger2}, and the present article is another illustration, even though we do quite a bit of number-crunching as well.

%%%%%%%%%%%%%%%%%%%%%%%%%%%%%%%%%%%%%%%%%%%%%%%%%%%%%%%%%%%%%%%%%%%%%%%%%%%%%%%%%%%%%%%%%%%%%%%%%
% Section: The Backward Induction Algorithm
%%%%%%%%%%%%%%%%%%%%%%%%%%%%%%%%%%%%%%%%%%%%%%%%%%%%%%%%%%%%%%%%%%%%%%%%%%%%%%%%%%%%%%%%%%%%%%%%%
\section{ The Backward Induction Algorithm}

The reason that it is so hard to decide (in some cases, for example with 5 heads and 3 tails) whether to {\bf stop} (and collect, {\it for sure}, the current number of heads divided by the current number of tosses [i.e. $h/(h+t)$]), or to  keep {\bf go}ing, ({\it expecting} to do better), 
is the somewhat unrealistic {\it assumption} that {\bf we live for ever}. Since in {\it real} life, we eventually would have to quit playing after $N$ tosses, for some finite $N$, and collect whatever we get then. So let's consider the {\it bounded} case where the number of coin-tosses is $\leq N$, for a fixed, possibly large, yet finite $N$. Compromising however with our immortality fantasy, we will let the player collect $1/2$, once reaching
the $N$-th coin toss, if the number of tails exceeds the number of heads, citing the law of large numbers that ``guarantees''
that ``eventually'' we will be able to (at least) break even. In other words, we let people who die in debt take advantage of the law of large numbers down in hell. [It turns out that, as far as the 
soon-to-be-defined limit, $F(h,t)$ goes, one does not need this assumption, and it is possible to insist that the player 
collects $h/N$ no matter what, but the breaking-even assumption considerably accelerates the convergence.]

Let's call $f_N(h,t)$ the expected pay-off in this bounded game, if you currently have $h$ heads and $t$ tails. Following Chow and Robbins, there is a simple backward induction (dynamical programming) algorithm for computing $f_N(h,t)$ 
for {\it all} $(h,t)$ with $h+t \leq N$.\\

\noindent
{\bf Boundary conditions}: when $h+t=N$:
\begin{equation}
f_N(h,N-h)=\max(1/2,h/N) \quad , \quad  (0 \leq h \leq N) .
\end{equation}

\noindent
{\bf Backward Induction}:
\begin{equation}
f_N(h,t)= \max \, \left (\, \frac{f_N(h+1,t)+f_N(h,t+1)}{2}, \quad \frac{h}{h+t} \, \right ).
\end{equation}

\noindent
[If you keep going, the expected gain is $[f_N(h+1,t)+f_N(h,t+1) ]/2$, if you stop the expected (and actual) gain is $h/(h+t)$. ]\\

\noindent
[$f_N(h,t)$ is implemented in procedure {\tt CR(h,t,N) } in {\tt ChowRobbins}. {\tt CRm(h,t,N) } is a faster version].\\

It is obvious that, for each specific $h$ and $t$, $f_N(h,t)$ is an increasing sequence in $N$, bounded above by $1$, so we know that the limit
\begin{equation}
F(h,t):=\lim_{N \rightarrow \infty} f_N(h,t) \quad ,
\end{equation}
``exists''.

Fantasizing that we {\it actually know} the values of $F(h,t)$,  (as opposed to knowing that they ``exist''), we can decide whether to {\bf stop} or {\bf go}. If $F(h,t)=h/(h+t)$ then we {\bf stop}, and otherwise we {\bf go}. This assumes that the player merely evaluates situations by {\it expectation}. As we know from the St. Petersburg paradox, expectation is not everything, and a player may choose to  {\it guarantee} collecting $h/(h+t)$ rather than taking a huge chance of eventually getting less. We will later describe other criteria for stopping.

Julian Wiseman \cite{wiseman1} estimates $F(0,0)$  to be $0.79295350640 \dots$ .

The difficulty in {\it proving}, for a given number of heads and tails, $(h,t)$, that it is optimal to stop is that we need {\it rigorous} non-trivial (i.e. $<1$) upper bounds valid for $f_N(h,t)$ for all $N$. Then this would also be true of $F(h,t)$, the limit as $N \rightarrow \infty$ of $f_N(h,t)$.
On the other hand it is easy to come up with lower bounds, namely $f_{N_0}(h,t)$ is $ \leq f_N(h,t)$ for all $N \geq N_0$, so in particular every specific $f_{N_0}(h,t)$ serves as a lower bound of $F(h,t)$, so it follows that whenever, for some $N_0$, it is true that $h/(h+t)<f_{N_0}(h,t)$, then we know {\it for sure} that it is good to {\bf go}.

%%%%%%%%%%%%%%%%%%%%%%%%%%%%%%%%%%%%%%%%%%%%%%%%%%%%%%%%%%%%%%%%%%%%%%%%%%%%%%%%%%%%%%%%%%%%%%%%%
% Section: The (probable) sequence
%%%%%%%%%%%%%%%%%%%%%%%%%%%%%%%%%%%%%%%%%%%%%%%%%%%%%%%%%%%%%%%%%%%%%%%%%%%%%%%%%%%%%%%%%%%%%%%%%
\section{The (probable) sequence $\beta_n$ }

So let's be realistic and take $N$ to be $50000$, rather than $\infty$. The sequence $\beta_n(50000)$, that we conjecture equals the ``real thing'' $\beta_n=\beta_n(\infty)$, for $1 \leq n \leq 185$, equals:
$$
1, 2, 3, 2, 3, 2, 3, 2, 3, 4, 3, 4, 3, 4, 3, 4, 5, 4, 5, 4, 5, 4, 5, 4, 5, 4, 5, 4, 5, 6, 5, 6, 5, 6, 5, 6, 5, 6,
$$
$$
5, 6, 5, 6, 7, 6, 7, 6, 7, 6, 7, 6, 7, 6, 7, 6, 7, 6, 7, 6, 7, 8, 7, 8, 7, 8, 7, 8, 7, 8, 7, 8, 7, 8, 7, 8, 7, 8,
$$
$$
7, 8, 9, 8, 9, 8, 9, 8, 9, 8, 9, 8, 9, 8, 9, 8, 9, 8, 9, 8, 9, 8, 9, 8, 9, 10, 9, 10, 9, 10, 9, 10, 9, 10, 9, 10,
$$
$$
9, 10, 9, 10, 9, 10, 9, 10, 9, 10, 9,10, 9, 10, 9, 10,11, 10, 11, 10, 11, 10, 11, 10, 11, 10, 11, 10,
$$
$$ 
11, 10, 11, 10, 11, 10, 11, 10, 11, 10, 11, 10, 11, 10, 11, 12, 11, 12, 11, 12, 11, 12, 11, 12, 11, 12,
$$
$11, 12, 11, 12, 11, 12, 11, 12, 11, 12, 11, 12, 11, 12, 11, 12, 11, 12, 11$. \\\\
We observe that for $1 \leq n \leq 9$, $\beta_{n^2}=n$ while for $10 \leq n \leq 13$, it equals $n-2$. This seems to be in harmony with Shepp's
theorem,  even for small $n$.

%%%%%%%%%%%%%%%%%%%%%%%%%%%%%%%%%%%%%%%%%%%%%%%%%%%%%%%%%%%%%%%%%%%%%%%%%%%%%%%%%%%%%%%%%%%%%%%%%
% Section: The question of when to stop and when to go depends on how long you expect to live.
%%%%%%%%%%%%%%%%%%%%%%%%%%%%%%%%%%%%%%%%%%%%%%%%%%%%%%%%%%%%%%%%%%%%%%%%%%%%%%%%%%%%%%%%%%%%%%%%%
\section{The question of when to stop and when to go depends on how long you expect to live}

We mentioned above that Ted Hill \cite{hill1} erroneously stated that 2 heads and 1 tails is a {\bf stop}. Well, he was not completely wrong. With $N \leq 50$, in other words, if the game lasts at most $50$ rounds, and as soon as you have tossed the coin $50$ times you must collect $\max(1/2,h/50)$, then $(2,1)$ is indeed a {\bf stop}. However, if the duration of the game is $\geq 51$, then it becomes a {\bf go}. We say that the {\it cutoff} for $(2,1)$ is $51$. In the following list, the $i$-th item is a pair. Its first component is
that position with $h+t=i$ that has 
the largest $h$ for which $(h,t)$ 
is a {\bf go} (for $N=2000$, and most probably (but unprovably) for $N=\infty$). Its second component is the smallest $N$ for which it stops being {\bf stop} and starts being {\bf go}. Notice the cautionary tales of the position with 10 heads and 7 tails that only starts being a {\bf go} with $N=1421$, and the position with 24 heads and 19 tails, for which $N=1679$ is the start of {\bf go}-dom.

Here is the list of pairs:

$$
[[[0, 1], 2], [[1, 1], 3], [[2, 1], 51], [[2, 2], 5], [[3, 2], 7], [[3, 3], 7],
    [[4, 3], 9], [[4, 4], 9], [[5, 4], 11],
$$
$$
    [[6, 4], 35], [[6, 5], 13], [[7, 5], 23], [[7, 6], 15], [[8, 6], 21], [[8, 7], 17], [[9, 7], 21],
    [[10, 7], 1421],
$$
$$
    [[10, 8], 23], [[11, 8], 91], [[11, 9], 25], [[12, 9], 57], [[12, 10], 25], [[13, 10], 47], [[13, 11], 27],
$$
$$
    [[14, 11], 43], [[14, 12], 29], [[15, 12], 43], [[15, 13], 31], [[16, 13], 43], [[17, 13], 277], [[17, 14], 43], 
$$
$$
    [[18, 14], 139], [[18, 15], 43],
    [[19, 15], 103], [[19, 16], 45], [[20, 16], 87], [[20, 17], 45], [[21, 17], 79],
$$
$$
    [[21, 18], 47], [[22, 18], 75], [[22, 19], 49],
    [[23, 19], 73], [[24, 19], 1679], [[24, 20], 71], [[25, 20], 423],
$$
$$
    [[25, 21], 71], [[26, 21], 249], [[26, 22], 69], [[27, 22], 185],
    [[27, 23], 69], [[28, 23], 155], [[28, 24], 71],
$$
$$
    [[29, 24], 137], [[29, 25], 71], [[30, 25], 125], [[30, 26], 73], [[31, 26], 119],
    [[31, 27], 73], [[32, 27], 113], 
$$
$$
    [[32, 28], 75], [[33, 28], 109], [[34, 28], 833], [[34, 29], 107], [[35, 29], 477], [[35, 30], 107],
    [[36, 30], 343],
$$
$$
    [[36, 31], 105], [[37, 31], 275], [[37, 32], 105], [[38, 32], 235], [[38, 33], 105], [[39, 33], 211], [[39, 34], 105],
$$
$$
    [[40, 34], 193], [[40, 35], 105], [[41, 35], 181], [[41, 36], 105], [[42, 36], 171], [[42, 37], 105], [[43, 37], 165],
$$
$$
    [[43, 38], 107],
    [[44, 38], 159], [[45, 38], 1039], [[45, 39], 155], [[46, 39], 679], [[46, 40], 153], 
$$
$$
[[47, 40], 513], [[47, 41], 151], [[48, 41], 419],
    [[48, 42], 149], [[49, 42], 361], [[49, 43], 147], [[50, 43], 321], 
$$
$[[50, 44], 147], [[51, 44], 293], [[51, 45], 147], [[52, 45], 271],[[52, 46], 145], [[53, 46], 255],$\\
$[[53, 47], 145] ].$
%%%%%%%%%%%%%%%%%%%%%%%%%%%%%%%%%%%%%%%%%%%%%%%%%%%%%%%%%%%%%%%%%%%%%%%%%%%%%%%%%%%%%%%%%%%%%%%%%
% Section: More Statistical Information
%%%%%%%%%%%%%%%%%%%%%%%%%%%%%%%%%%%%%%%%%%%%%%%%%%%%%%%%%%%%%%%%%%%%%%%%%%%%%%%%%%%%%%%%%%%%%%%%%
\section{More Statistical Information}

The above strategy for deciding when to stop  is entirely based on {\it expectation}. 
Even if we pursue this strategy, it would be nice to have more detailed information, like the {\it standard deviation}, {\it skewness}, {\it kurtosis} and even higher moments. Ideally, we would like to know the {\it full probability distribution}. 

Let's call $G_N(h,t;x)$ the fractional polynomial in the variable $x$ (i.e. a linear combination of powers $x^a$ with $a$ rational numbers) such that the coeff. of $x^a$ is the probability of getting exactly $a$ as pay-off in our game, still pursuing the strategy of maximizing the expected gain.  Of course $G_N(h,t;1)=1$ and ${{d} \over {dx}} G_N(h,t;x) \vert_{x=1}=f_N(h,t)$. We have:\\

\noindent
{\bf Boundary conditions}: when $h+t=N$:
\begin{equation}
G_N(h,N-h;x)=x^{\max(1/2,h/N)} \quad  (0 \leq h \leq N) \quad.
\end{equation}
{\bf Backward Induction}:

\begin{eqnarray}
\quad
G_N(h,t;x)=
\begin{cases}
x^{h/(h+t)} ,& \text{if ({\it h},{\it t}) is STOP} \\
{{G_N(h+1,t;x)+G_N(h,t+1;x)} \over {2}},& \text{if ({\it h},{\it t}) is GO.}
\end{cases},
\end{eqnarray}

\noindent
[$G_N(h,t;x)$ is implemented in procedure {\tt CRt(h,t,N,x) } in {\tt ChowRobbins}.]\\

\noindent
Once we have $G_N(h,t;x)$, we can easily get all the desired statistical information.

%%%%%%%%%%%%%%%%%%%%%%%%%%%%%%%%%%%%%%%%%%%%%%%%%%%%%%%%%%%%%%%%%%%%%%%%%%%%%%%%%%%%%%%%%%%%%%%%%
% Section: Another Way to Gamble
%%%%%%%%%%%%%%%%%%%%%%%%%%%%%%%%%%%%%%%%%%%%%%%%%%%%%%%%%%%%%%%%%%%%%%%%%%%%%%%%%%%%%%%%%%%%%%%%%
\section{Another Way to Gamble}

In real life we don't always want to maximize our expected gain. Often we have a certain {\it goal}, let's call it $g$, and achieving or exceeding it means everlasting happiness, while getting something less would mean eternal misery. In that case we need a different gambling strategy, that is really straightforward. Keep playing until $h/(h+t) \geq g$, and if and when you reach it, {\bf stop}. Otherwise keep going to the end, until $h+t=N$. In that case, of course, the {\bf stop} states are those for which $h/(h+t) \geq g$. It is still of interest to to know what {\it is} the probability of happiness. 
Let's call this quantity $P_N(g;h,t)$. We obviously have:\\

\noindent
{\bf Boundary conditions}: when $h+t=N$:
\begin{eqnarray}
P_N(g;h,N-h)=\begin{cases}
0 ,& \text{if } h/N < g\\
1,&  \text{if } h/N \geq  g.
\end{cases},
\end{eqnarray}

\noindent
{\bf Backward Induction}: When $h+t<N$,
$P_N(g;h,t)$ equals $1$ if $h/(h+t) \geq g$ while it equals
$(P_N(g;h+1,t)+P_N(g;h,t+1))/2$ otherwise.\\

We leave it to the reader to formulate the backward induction scheme for finding
the probability generating function for the present strategy.

%%%%%%%%%%%%%%%%%%%%%%%%%%%%%%%%%%%%%%%%%%%%%%%%%%%%%%%%%%%%%%%%%%%%%%%%%%%%%%%%%%%%%%%%%%%%%%%%%
% Section: Comparative Gambling
%%%%%%%%%%%%%%%%%%%%%%%%%%%%%%%%%%%%%%%%%%%%%%%%%%%%%%%%%%%%%%%%%%%%%%%%%%%%%%%%%%%%%%%%%%%%%%%%%
\section{Comparative Gambling}

Let's  compare the two strategies using {\it both} criteria. Of course the first one always is better in the {\it maximum expectation category} and the second is always better in {\it maximizing the probability of achieving the goal}.

With $N=200$, at the very beginning, your expected gain, under the first way is  $0.7916879464$, but your probability
\begin{itemize}
\item of getting $\geq 0.6$ is $0.6917238235$ (the second way gives you probability $0.7753928313$, but your expected gain is only $0.6742902054$)

\item of getting $\geq 0.7$ is $0.5625000000$ (the second way gives you probability $0.6075176458$, but your expected gain is only $0.5787939263$)
\end{itemize}

\noindent
Much more data can be found by using procedure {\tt SipurCG} in the Maple package {\tt  ChowRobbins}, and posted in the webpage of this article.

%%%%%%%%%%%%%%%%%%%%%%%%%%%%%%%%%%%%%%%%%%%%%%%%%%%%%%%%%%%%%%%%%%%%%%%%%%%%%%%%%%%%%%%%%%%%%%%%%
% Section: Probabilities of Escape
%%%%%%%%%%%%%%%%%%%%%%%%%%%%%%%%%%%%%%%%%%%%%%%%%%%%%%%%%%%%%%%%%%%%%%%%%%%%%%%%%%%%%%%%%%%%%%%%%
\section{Probabilities of Escape}

The second strategy gives rise to the following interesting computational question:\\

\noindent
Fix $a>b \geq 1$ relatively prime. What is the probability that the number of heads divided by the number of tails

(i) will ever exceed $a/b$? 

(ii) will either exceed or be equal to $a/b$?\\

\noindent
This question was raised and answered by Wolfgang Statdje \cite{stadje1} who proved that this quantity is a root of a certain algebraic equation. A related problem is treated by Nadeau \cite{nadeau1}.

Stadje's result can also be deduced from the more general treatment by Ayyer and Zeilberger \cite{AZ}, that contains a Maple package that automatically derives the algebraic equation for any general set of steps. For practical purposes, however, we found it easiest to compute these
probabilities directly, in terms of the discrete functions $W(x,y)$ and $W_s(x,y)$ that count the number of lattice walks from the origin to $(x,y)$ staying in the required region. This is contained in the Maple package {\tt STADJE}.

Here is some data gotten from {\tt STADJE}. The numbers below answer questions (i) and (ii) above, respectively, for each of the listed pairs $(a,b)$.\\

$(a,b)=(2,1):  0.6180339887, 0.6909830056$ ;

$(a,b)=(3,1):  0.5436890127, 0.5803566224$ ;

$(a,b)=(3,2):  0.7481518342, 0.7754441182$;

$(a,b)=(4,1):  0.5187900637,  0.5362190123$ ;

$(a,b)=(4,3):  0.8091410707, 0.8229424412 $;

$(a,b)=(5,1):  0.5086603916, 0.5170258817 $;

$(a,b)=(5,2):  0.5876238826, 0.5996923731$;

$(a,b)=(5,3):  0.7158769909, 0.7276461121 $;

$(a,b)=(5,4):  0.8453136528, 0.8534748833 $;\\

Also of interest is the sequence enumerating the number of walks, staying in the region $y \geq a/bx$,  from the origin to a point of the form $(n,n)$, whose asymptotics can be proved to be of the form $C_1(a,b)4^n/\sqrt{n}$, for some constant $C_1(a,b)$, and the sequence enumerating the number of walks,
still staying in the same region, ending at $(an,bn)$, whose asymptotics has the form $C_2(a,b) ((a+b)^{a+b}/(a^ab^b))^n/n^{3/2}$. 
The Maple package {\tt STADJE} (and Mathematica package {\tt STADJE.m})
computes any desired number of terms, and estimates $C_1(a,b)$, $C_2(a,b)$. 
The webpage of this article contains some sample output.

%%%%%%%%%%%%%%%%%%%%%%%%%%%%%%%%%%%%%%%%%%%%%%%%%%%%%%%%%%%%%%%%%%%%%%%%%%%%%%%%%%%%%%%%%%%%%%%%%
% Section: From Number-Crunching to Symbol Crunching
%%%%%%%%%%%%%%%%%%%%%%%%%%%%%%%%%%%%%%%%%%%%%%%%%%%%%%%%%%%%%%%%%%%%%%%%%%%%%%%%%%%%%%%%%%%%%%%%%
\section{From Number-Crunching to Symbol Crunching}

So far, we have designed {\it numerical} computer programs whose outputs were {\it numbers}. But what about {\it closed form}? It would be too much to hope for an explicit formula for $f_N(h,t)$ valid for {\it arbitrary} $N$, $h$, $t$, {\it but}, with {\it experimental-yet-rigorous} mathematics, we can find {\it explicit expressions}, as {\it rational functions} in $n$ for
\begin{equation}
f_{2n+1}(n+\alpha,n-\alpha-m+1),
\end{equation}
where $n$ and $m$ are positive integers and $\alpha$ is an integer.
%Of course when $\alpha=0$, $F_N(h,N-h)=\max(1/2,h/N)$, by definition, but what about larger $\alpha$?

%Our focus in this section is to get closed formulas for $f(h,t,n)$ when $n$ is big.  In particular, we study 
Let
\begin{equation}
F(m,\alpha,n)=f_{2n+1}(n+\alpha,n-\alpha-m+1)
\end{equation}
for $n$, $m$, and $\alpha$ as before.  Since $h+t < 2n+1$, then $F(m,\alpha,n)$ are values below the topmost diagonal on the backward induction triangle.

Some values of $F(m,\alpha,n)$ are not hard to get. For instance, the value of $F(m,\alpha,n)$, for $\alpha\geq 1$ and $1\leq m\leq2n$, is given by
\begin{equation}
F(m,\alpha,n)=\frac{n+\alpha }{2n-m+1},
\end{equation}
whereas the value of $F(m,\alpha,n)$, for $\alpha\leq -m$ and $1\leq m \leq 2n$, is given by
\begin{equation}
F(m,\alpha,n)=\frac{1}{2}.
\end{equation}
Both formulas can be proved by induction.  Hence, we are reduced to finding formulas for $F(m,\alpha,n)$ when $-m < \alpha <1.$

Our first approach is to make the computer conjecture closed forms for $F(m,\alpha,n)$.  For this, we programmed a Mathematica function called \texttt{GF} [this function can be found in the webpage of this article].  It takes as input a positive integer $m$ and two variables $n$ and $\alpha$, and another positive integer $bound$.  Here, the computer makes the assumption that $n\geq bound$.  For the guessing part, \texttt{GF} uses the auxiliary function \texttt{GuessRationalFunction}.  This procedure is similar to \texttt{GuessRat}, which accompanied the article \cite{goodWay} and can be found in \cite{DZprograms}.  The output of \texttt{GF}, which is the guess formula for $F(m,\alpha,n)$,  is a piecewise rational function of $n$ with $m+2$ pieces.

\begin{example}
For $m=2$ and $n\geq 3$, \texttt{GF} conjectures
\begin{equation}
F(2,\alpha,n)=\begin{cases}
 1/2 & \alpha \leq -2 \\
 \frac{8 n+5}{16 n+8} & \alpha =-1 \\
 \frac{8 n^2+9 n+2}{16 n^2+8 n} & \alpha =0 \\
 \frac{n+\alpha }{2 n-1} & \alpha \geq 1
\end{cases}
\end{equation}
\end{example}

We point out that formulas conjectured by \texttt{GF} only work for $n$ sufficiently large.  In fact, empirical evidence suggests that the bound on $n$ grows exponentially in $m$ i.e. as we go down on the backward induction triangle, the bound for which the formulas are valid grows exponentially.  
As a result, these formulas are not {\it directly} useful for determining {\bf stop} vs. {\bf go} status.

%\section{Computer generated formulas and proofs}
It is possible to study the recursion formula of $f_n(h,t)$ to get explicit formulas for $F(m,\alpha,2n+1)$.  For example, a simple analysis gives
\begin{equation}
F(1,\alpha,2n+1)=\begin{cases}
 1/2 & \alpha \leq -1 \\
 \frac{4 n+3}{8 n+4} & \alpha =0 \\
 \frac{n+\alpha }{2 n} & \alpha \geq 1
\end{cases}
\end{equation}
which is true for $n\geq 1$, and
\begin{equation}
F(2,\alpha,2n+1)=\begin{cases}
 1/2 & \alpha \leq -2 \\
 \frac{8 n+5}{16 n+8} & \alpha =-1 \\
 \frac{8 n^2+9 n+2}{16 n^2+8 n} & \alpha =0 \\
 \frac{n+\alpha }{2 n-1} & \alpha \geq 1
\end{cases}
\end{equation}
which is true for $n\geq 3$.  However, these calculations become tedious rapidly.

To our surprise, it turns out that Mathematica, via the built-in functions \texttt{Assuming} and \texttt{Refine}, is able to handle these recursions and get the desired formulas.  We programmed a Mathematica function called \texttt{BUILDER}, whose input is an integer $m$ and two variables $n$ and $\alpha$.  \texttt{BUILDER} calculates closed-form formulas for $F(m,n,\alpha)$ and provides the smallest $n$ where they start to hold.  For instance,
\begin{equation}
F(5,\alpha,2n+1)=\begin{cases}
 1/2 & \alpha \leq -5 \\
 \frac{64 n+33}{128 n+64} & \alpha =-4 \\
 \frac{32 n^2+20 n+1}{64 n^2+32 n} & \alpha =-3 \\
 \frac{64 n^3+30 n^2-13 n-3}{128 n^3-32 n} & \alpha =-2 \\
 \frac{64 n^4+8 n^3-46 n^2-5 n+3}{128 n^4-128 n^3-32 n^2+32 n} & \alpha =-1 \\
 \frac{256 n^5-124 n^4-340 n^3+91 n^2+75 n-6}{512 n^5-1280 n^4+640 n^3+320 n^2-192 n} & \alpha =0 \\
 \frac{n+\alpha }{2 n-4} & \alpha \geq 1
\end{cases}
\end{equation}
was calculated by \texttt{BUILDER} and holds for $n \geq 102$.

The starting places, for $n$, where
the formulas of $F(m,\alpha,n)$ begin to hold,
with $1\leq m \leq 16$, are:  1, 3, 12, 37, 102, 263, 648, 1545, 3594, 8203, 18444, 40973, 90126, 196623, 426000, and 917521 respectively.  These values seems to satisfy the recurrence defined by
\begin{eqnarray*}
a_1 &=& 1 \\
a_m &=& 2 a_{m-1}+r_m \text{ valid for }m\geq 1,
\end{eqnarray*}
where $r_m$ is given by
\begin{eqnarray*}
r_1 &=& 0\\
r_2 &=& 1\\
r_3 &=& 6 \\
r_m &=& 2 r_{m - 1} + m - 3 \text{ valid for }m\geq3.
\end{eqnarray*}

We are pleased to report that the formulas conjectured by \texttt{GF} and the ones found by \texttt{BUILDER} agree.

%%%%%%%%%%%%%%%%%%%%%%%%%%%%%%%%%%%%%%%%%%%%%%%%%%%%%%%%%%%%%%%%%%%%%%%%%%%%%%%%%%%%%%%%%%%%%%%%%
% Acknowledgments
%%%%%%%%%%%%%%%%%%%%%%%%%%%%%%%%%%%%%%%%%%%%%%%%%%%%%%%%%%%%%%%%%%%%%%%%%%%%%%%%%%%%%%%%%%%%%%%%%
\bigskip

\noindent
{\bf Acknowledgments}. We wish to thank Theodore Hill for very useful Email correspondence.
The work of the second author was supported in part by the USA National Science Foundation.
\bigskip

\noindent
{\bf July 23, 2009}
\bigskip

%%%%%%%%%%%%%%%%%%%%%%%%%%%%%%%%%%%%%%%%%%%%%%%%%%%%%%%%%%%%%%%%%%%%%%%%%%%%%%%%%%%%%%%%%%%%%%%%%
% References
%%%%%%%%%%%%%%%%%%%%%%%%%%%%%%%%%%%%%%%%%%%%%%%%%%%%%%%%%%%%%%%%%%%%%%%%%%%%%%%%%%%%%%%%%%%%%%%%%

\end{document}